\documentclass[leqno,fleqn,a4paper]{amsart}

\usepackage{destemacros}

\begin{document}

\title[Recent and less recent results on Tilting Theory]{A personal survey on recent and less recent 
                             results on Tilting Theory}


\author{Gabriella D'Este}
\address{Department of Mathematics, University of Milano\\
	Via Saldini 50, 20133 Milano, Italy
}
\email{gabriella.deste@unimi.it}

\dedicatory{Dedicated to the memory of Sheila Brenner and  Michael C.~R.~Butler}


\date{}


\maketitle

\section*{Introduction}
In the introduction to Ringel's book ``Tame Algebras 
and Integral Quadratic Forms'', the Author makes 
the following comment about references: 
``however we should point out that some general 
ideas, which have influenced the results and the 
methods presented here, are not available in official 
publications, or not even written up.''  
The above remark also explains very well my personal 
experiences with respect to many situations concerning 
both Places and People, to use the terminology of  
FDLIST \cite{F} .  For me one of the first examples of the 
complexity of ``general ideas'' (and their zigzag journeys) 
is given by the beginning of Tilting Theory in Italy.  
The best written reference for the knowledge of the 
big role played by Adalberto Orsatti  - and his Algebra 
Team in Padova - is  Menini's paper \cite{M} .  This paper 
contains a very interesting account of both important 
public events and official publications, as well as 
information, never previously written up, on private 
conversations, classical letters and unexpected 
connections between distant places and people. On the 
other hand, the first paper by Italian authors on Tilting 
Theory, that is \cite{MO}, is due to Menini and Orsatti.  
The title of the paper (``Representable equivalences 
between categories of modules and applications''), 
suggests that classes of modules are the most 
important ingredients. In the introduction of  [MO]  
the authors thank Masahisa Sato, Enrico Gregorio and me.  
I regret that I have never expressed officially - in a paper - 
my thanks (and surprise) for this unexpected reference, as 
I should have done long ago. The aim of this paper is to fill 
this gap, at least partially.  In a sense, this note is the 
written version of conversations with young colleagues on 
unofficial history, ``general ideas'', unexpected facts  and 
open problems.  I may sum up as follows the lessons 
learned by making pictures of tilting--type objects:
\begin{enumerate}[(a)]
 \item ``Simple'' and combinatorial objects may have 
unexpected concealed topological properties.
\item ``Non simple'' objects may have unexpected 
concealed discrete properties.
\end{enumerate}
   
In the following  $K$  denotes an algebraically closed field,
and we assume that all vector spaces and algebras are 
defined over  $K$ .  This note is organized as follows.  In 
Section 1 , I will describe my first homework on Tilting 
Theory in a non technical way.  This homework, given 
by Orsatti, was to study and explain the example of 
tilting module constructed by Happel and Ringel at the
end of  \cite{HR}, hence ``coming from Bielefeld'', as Orsatti 
told me.  Next, in Sections 2  and  3, I describe the 
na\"\i{}ve strategies, used, from the very beginning, to 
answer some questions on more or less abstract tilting 
objects, and I make some personal remarks on definitions.  
Finally, in Sections  4  and  5, I collect some motivations 
and examples for the observations and/or conjectures 
made in (a)  and (b).
In addition to the names of  C. Menini  and A. Orsatti 
(the authors of  \cite{MO}  already mentioned), any partial 
list of people I would like to thank should contain at 
least the following names: L. Angeleri - H\"ugel, S. Bazzoni, 
R. Colpi, E. Gregorio, F. Mantese, A. Tonolo and J. Trlifaj.  
Thanks to all of you for your problems -- explained in the 
clearest possible way -- and your hope that there ought to 
be a solution somewhere!


I described some facts and results contained in this note in some 
seminars at the universities of Milano (June 2012), Milano Bicocca 
(October 2012) and Ferrara (June 2013). I wish to thank the 
organizers of these talks (Gilberto Bini, Cecilia Cavaterra, Elisabetta 
Rocca, Maria Gabriella Kuhn and Claudia Menini) for the opportunity 
of presenting my work to an interested audience, and all the 
participants for their interest and questions. Next, I would like to 
thank Ibrahim Assem, together with all the organizers of the 
"XVIIth Meeting on representation theory of algebras" (Bishop's 
University, Sherbrooke, October 2013) dedicated to the memory 
of Michael C.R. Butler, for the idea of putting a preliminary version
of this note in the home page of the meeting, with the title
 "Preprint dedicated to Brenner Butler". I will always remember with
gratitude many conversations with the creators of tilting theory 
and their attention to the work of other people of any age and 
countries.

\section{A letter from Japan and my first homework 
      on Tilting Theory}
Menini's paper \cite[page  11]{M} explains the contribution 
given by Sato to her joint paper  \cite{MO}  with Orsatti in
the following comments:
``I woud like to recall here that it was Masahisa Sato 
that pointed out to us that tilting modules might 
provide examples of \dots.    In fact  Orsatti explained 
this problem to Sato during a NATO meeting É held 
in Antwerpen (Belgium) in the period  July 20 - 29 ,
1987 .  After some time Sato wrote to Orsatti showing
an example of a tilting module''.
I remember very well what happened next.  Orsatti 
told me to make a copy of the letter received from 
Japan.  My homework was to look at Happel - Ringel's 
example of a tilting module, say  $T$, considered 
in the last page of \cite{HR}, and to give a talk about it. 
Hence my unique and small contribution to  \cite{MO}  
was just a talk.  This explains my surprise for those
unexpected thanks.   Moreover, now I realize that I 
should have thanked Menini and Orsatti.  Indeed the 
preparation of the seminar was very useful for me 
for several reasons.  For instance, I had the pleasant
surprise of a direct experience that  Auslander--Reiten 
quivers can really help to guess and see possible 
equivalences, before making a proof of their existence.  
Indeed tilting equivalences and cotilting dualities  
also have a combinatorial nature, inherited and suggested by 
that of quivers and modules.  In addition to this, I was
able to discover the magic power of  Auslander's formula 
(see  \cite[Proposition 4.6  and Corollary 4.7]{AuReS}   or  
\cite[conditions  (5)  and (6) , pages  75 - 76]{R4})  for verifying
the vanishing of certain  $\Ext^1$      groups, by simply looking at 
the Auslander--Reiten quiver.  Without this formula, I am
(and was) not able to check that the module  $T$  is
selforthogonal.  The next picture (of a very combinatorial 
object with several symmetries) illustrates the shape of  
Happel--Ringel's bimodule  $_AT_B $.
\eject
\begin{center}
\textbf{\large Picture 1} \\
\bigskip 
\hspace*{-20mm}
\rotatebox{90}{
$$\scalebox{0.95}{$
\xymatrix{
						& & & & & \fquadrcolor{5}{a}{green} \ar@{-}[d] \ar@{~}@/_2pc/[ddddd] \ar@{~}[ddrrr] \ar@{~}[ddlll] &  \\
						& & & & & \fquadrcolor{6}{a}{cyan} \ar@{~}[ddddddr]\ar@{~}[ddddddl] \ar@{~}[ddllll]\ar@{~}[ddrrrr] &  \\
\fquadrcolor{2}{b}{yellow} \ar@{-}[dr]\ar@{~}@/_3pc/[ddddddd]& & \fquadrcolor{5}{b}{green} \ar@{~}[ddddddd] \ar@{-}[dl]& & & & & & \fquadrcolor{5}{c}{green} \ar@{-}[dr]\ar@{~}[ddddddd]  & & \fquadrcolor{4}{c}{orange} \ar@{-}[dl]\ar@{~}@/^3pc/[ddddddd] \\
	&  \fquadrcolor{6}{b}{cyan} \ar@{~}@/_3pc/[ddddddd]  &  & & & & & &  & \fquadrcolor{6}{c}{cyan}\ar@{~}@/^3pc/[ddddddd] & \\
	& & & \fquadrcolor{1}{d}{red} \ar@{~}[ddddlll] \ar@{-}[dd] & & & & \fquadrcolor{3}{d}{magenta}  \ar@{~}[ddddrrr] \ar@{-}[dd] & & & \\
	& & &        & & \fquadrcolor{5}{d}{green} \ar@{~}@/_1pc/[ddddlll]\ar@{~}@/^1pc/[ddddrrr] \ar@{-}[ddr]\ar@{-}[ddl]& & & & & \\
	& & & \fquadrcolor{2}{d}{yellow} \ar@{~}[dddlll] \ar@{-}[dr] & & & & \fquadrcolor{4}{d}{orange} \ar@{~}[dddrrr] \ar@{-}[dl]& & & \\
	& & & & \fquadrcolor{6}{d}{cyan}  \ar@{~}@/^3pc/[dddlll] & & \fquadrcolor{6}{d}{cyan} \ar@{~}@/_3pc/[dddrrr] & & & & \\
	\fquadrcolor{1}{e}{red} \ar@{-}[d]& & & & & & & & {\phantom{\fquadrcolor{1}{e}{red}}} & & \fquadrcolor{3}{f}{magenta} \ar@{-}[d]\\
	\fquadrcolor{2}{e}{yellow}  \ar@{-}[dr]&& \fquadrcolor{5}{e}{green} \ar@{-}[dl] &&&&&& \fquadrcolor{5}{f}{green}  \ar@{-}[dr]&& \fquadrcolor{4}{f}{orange} \ar@{-}[dl]\\
	& \fquadrcolor{6}{e}{cyan} & & & & & & & & \fquadrcolor{6}{f}{cyan} & 
\save "1,6"."2,6"*+<15pt>[F.:<10pt>]\frm{}\restore
\save "3,1"."4,3"*+<15pt>[F.:<10pt>]\frm{}\restore
\save "3,9"."4,11"*+<15pt>[F.:<10pt>]\frm{}\restore
\save "5,4"."8,8"*+<15pt>[F.:<10pt>]\frm{}\restore
\save "9,1"."11,3"*+<15pt>[F.:<10pt>]\frm{}\restore
\save "9,9"."11,11"*+<15pt>[F.:<10pt>]\frm{}\restore
}$}
$$
}
\end{center}

According to  \cite{HR} ,  the underlying vector space of  $T$  has
dimension  23 , while  $A$   is the algebra given by the Dynkin 
quiver  $E_6$     with ``subspace orientation'', that is of the form
            
$$
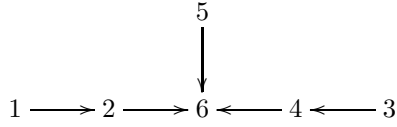
\xymatrix{
   &   & 5 \ar[d] & & \\
 1 \ar[r]& 2 \ar[r] & 6 & 4 \ar[l]& 3 \ar[l]
}            
$$
     
\bigskip

On the other hand, the algebra  $B  =  \End T_A $    is isomorphic 
to the algebra given by the fully commutative quiver

$$
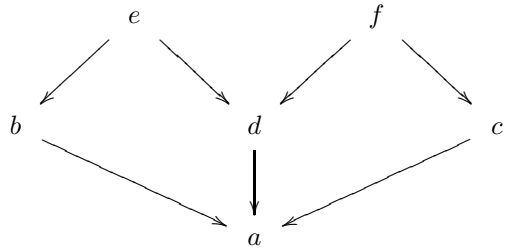
\xymatrix{
      & \coluno{e} \ar[dl]\ar[dr] & &\coluno{f} \ar[dl] \ar[dr] & \\
 \coluno{b} \ar[drr] &  &  \coluno{d} \ar[d]  &   & \coluno{c} \ar[dll] \\
 && \coluno{a} & &
}
$$

Finally, in Picture 1  and in the next pictures of bimodules, 
we adopt the following conventions.  First of all, every 
square of the picture indicates an element   $v$  of a fixed 
basis of the underlying vector space of  $T$ .  Next, the 
index   $x$  on the left  (resp.~ $y$  on the right)  of a small 
square corresponding to the vector  $v$  indicates that  
$e_x  v = v =  v e_y$   , where  $e_i$   is the path of length zero
around the vertex  $i$.  Following Ringel's suggestion
during my staying in Bielefeld, the small squares 
have a special position, so that they describe in an 
obvious way also the composition factors of the same      
module. (See, for instance, \cite{R1,R2}  and  \cite[page
126]{R3}  for descriptions and/or pictures of complicated 
modules.)  On the other hand, the bimodules associated 
to ``valued'' arrows in  \cite{DlR}  and  \cite{Dl}  gave me the idea 
of adding two indices (on the left and on the right of the 
small squares).  In this way we can see the action of 
left or right multiplication by the primitive idempotents, 
corresponding to vertices of some quiver.  Finally, I 
want to say that the above picture is not as old as 
my homework.  I used similar pictures to visualize 
rather small bimodules for the first time, and more 
or less by chance.  Of course, in case of big and 
complicated bimodules, it would be better to replace 
a  2--dimensional ``global'' visualization (of the main 
properties of the left and right underlying modules) 
by a  3--dimensional one, without or with few 
self - intersections. However, even flat pictures, 
like Picture 1 , are powerful enough to give for 
free a lot of  indirect information in a compact way.  
As an example, the above picture tells us how the 
tilting equivalence represented by  $T$ (between the 
modules generated by the left  $A$--module  $T$   and 
the modules cogenerated by the left  $B$--module  
$D(T_B) =  \Hom_K  ( T_B   , K )  $)     acts on some  
indecomposable modules.  For instance, the
indecomposable summands of   $_AT$   are sent 
to the following indecomposable summands of     $_BB$ :
                                                                  
$$
\coldue{5}{6}\mapsto a, \quad \coldue{2\ 5}{\ 6} \mapsto \coldue{b}{a}, \quad
\coldue{5\ 4}{\ 6}\mapsto \coldue{c}{a},
$$

$$
\coltre{1\bx\bx\bx\bx\bx3}{\bx2\bx5\bx4\bx}{\bx\bx6\bx6\bx\bx} \mapsto \coldue{d}{a}, \quad 
\coltre{1\bx\bx\bx}{\bx2\bx5}{\bx\bx6\bx} \mapsto \coltre{\bx{}e\bx }{b\bx{}d}{\bx{}a\bx}, \quad 
\coltre{\bx\bx\bx3}{5\bx4}{\bx6\bx} \mapsto \coltre{\bx{}f\bx}{d\bx{}c}{\bx{}a\bx{}}. 
$$  

On the other hand  Picture 1  also describes how the 
cotilting duality induced by  $T$ (between the modules 
cogenerated by left  $A$ - module  $T$    and the modules 
cogenerated by the right  $B$ - module  $T$ )  acts on some          
indecomposable modules.  For instance, the
indecomposable summands of    $_AT$     are sent
to following indecomposable summands of  $B_B$:
$$
\coldue{5}{6} \mapsto \coltre{\bx\bx{}a\bx\bx}{b\bx{}d\bx{}c}{\bx{}e\bx{}f\bx}, \quad 
\coldue{2\bx5}{\bx6\bx} \mapsto \coldue{b}{e}, \quad 
\coldue{5\bx4}{\bx6\bx} \mapsto \coldue{c}{f},
$$                                                                
$$
\coltre{1\bx\bx\bx\bx\bx3}{\bx2\bx5\bx4\bx}{\bx\bx6\bx6\bx\bx} \mapsto \coldue{\bx{}d\bx{}}{e\bx{}f}, \quad 
\coltre{1\bx\bx\bx}{\bx2\bx5}{\bx\bx6\bx} \mapsto \coluno{e}, \quad 
\coltre{\bx\bx\bx3}{5\bx4}{\bx6\bx} \mapsto \coluno{f}. 
$$

\section{What happened next: talk with big matrices, 
      more or less abstract cancellations\dots }
The analysis of  Happel--Ringel's bimodule was both 
my first homework on Tilting Theory, and the subject 
of my first talk on this subject.  The talk was not my 
first talk in Padova containing some quivers. However, 
all the quivers used earlier were much smaller.  
Hence, it was easier to describe without pictures at 
least  the corresponding path algebras.  When I tried 
to do the same with Happel--Ringel's example, I 
realized how even basic techniques of representation 
theory of finite dimensional algebras can make
otherwise invisible things become visible.  For instance, 
in order to give a definition of  A  and  B   without 
quivers, I distributed some pages with the largest 
matrices I ever used, namely  23  by  23  matrices, 
that I did not wanted to draw at the blackboard.  After
some time I noticed that I could use more reasonable 
matrices to describe the  $K$--linear maps from  $T$  to  
$T$, corresponding to multiplications by elements of
$A$  or  $B$ .  Indeed, $T$  has the strong property that 
the groups of  all morphisms between two 
indecomposable summands of   $_AT$   and   $T_B$  
respectively are  $K$--vector spaces of dimension at 
most one.  Hence, after cancellation of  many 
``inessential'' rows and columns, it is easy to see that  
$A$  and  $B$  are isomorphic to subalgebras of  the 
algebra of all  6  by  6  matrices of the following form:

$$
\begin{pmatrix} 
  K & 0 & 0 & 0 & 0 & 0  \\
  K & K & 0 & 0 & 0 & 0  \\
  0 & 0 & K & 0 & 0 & 0  \\
  0 & 0 & K & K & 0 & 0  \\
  0 & 0 & 0 & 0 & K & 0  \\
  K & K & K & K & K & K  \\
\end{pmatrix}
\qquad, \qquad
\begin{pmatrix} 
  K & K & K & K & K & K  \\
  0 & K & 0 & 0 & K & 0  \\
  0 & 0 & K & 0 & 0 & K  \\
  0 & 0 & 0 & K & K & K  \\
  0 & 0 & 0 & 0 & K & 0  \\
  0 & 0 & 0 & 0 & 0 & K  \\
\end{pmatrix}
.
$$

It turned out that endomorphism rings of abelian 
groups are the subject of my master degree thesis 
on ``Abelian groups whose endomorphism ring is 
locally compact in the finite topology'', written 
under the direction of  Adalberto Orsatti. Moreover, 
some of my first papers deal with endomorphism 
rings. However, without using quivers, this previous 
abstract experience on (usually large) endomorphism 
rings wouldn't have been useful for verifying that  
$B = \End _AT $   actually had the indicated form.  The  
cancellation of rows and columns was only the first
of many other (more abstract) cancellations made in 
the sequel, concerning modules and complexes, as in  
3.5 , 3.6  and  3.7.  For the theoretical importance 
of cancellations in Tilting Theory, we refer to Ringel's 
lecture on the occasion of the  20th anniversary of the 
Department of Mathematics of the University of Padova 
[see the section ``Fully documented Lectures'' in Ringel's 
home page].  We note that the suggestive title of this 
lecture is ``Tilting Theory: the Art of Losing Modules''.

\section{Short or long definitions with or without 
     classes of modules}
It is well--known that some definitions of tilting and/or 
partial tilting modules (the most compact and elegant 
ones) consist of precisely one rather short property.  
For instance, any ``classical'' tilting (resp.~cotilting) 
$R$--module  $M$  (hence also Happel--Ringel's module  $T$), 
has the property that the class of all modules 
generated (resp.~ cogenerated) by  $M$  coincides with 
the kernel of   $\Ext^1_R(M,-)$ (resp.~$\Ext^1_R(-,M)$).  This 
previous global property is equivalent to three discrete 
properties of two modules, namely of  $M$  and of the 
regular module  $R$  (resp.~an injective cogenerator  $Q$), 
as in the definition given by Brenner and Butler in \cite{BB} .
    
\subsection{A long definition (for the classical case).}  
We say that a module  $M$  is a classical tilting or cotilting
module (more precisely, a  1--tilting or  1--cotilting
module) respectively, if the following conditions hold:
\begin{itemize}
 \item The projective (resp.~injective) dimension of  $M$  is 
  at most  1.
  \item  $\Ext^1_R    ( M ,  \bigoplus  M ) = 0$  (resp.~ $\Ext^1_R   ( \Pi M , M ) = 0$ ), 
  where    $\bigoplus M$   (resp.~   $\Pi  M$) is any direct sum  
    (resp.~ product)  of copies of  $M$. 
    \item There is a short exact sequence of the form
  $0 \To  R \To M' \To M''\To 0$    (resp.
   $0 \To M'\To M''\To Q \To 0$) ,  where   
   $M'$ and $M''$  are direct summands of direct sums 
(resp.~ products) of copies of  $M$.
\end{itemize}
For me it was always easier to deal with more than one 
(but finitely many) elementary properties instead of dealing 
with just one property on classes of modules.   We refer to  
\cite[Proposition 3.6  and  Lemma 3.12]{Ba}   for the beautiful 
technical condition (on the relationship between two classes 
of modules) which characterizes ``non classical''  tilting or 
cotilting modules and their generalizations (for instance, their 
direct summands and more precisely partial tilting (resp.~
cotilting) modules of projective (resp.~injective) dimension  
$> 1$).  We only recall that - as in the classical case - the 
usual definition of all these modules  $M$  consists of two 
conditions on just two modules, namely on $M$  and a very 
special projective or injective module.  One of the reasons 
why it may be difficult to check equalities or inclusions of 
classes of modules which play a key role in the 
characterization of  ``non classical''  tilting and cotilting -
type modules  $M$  is the following.  Even in dealing with 
algebras of finite representation type, in the ``non classical'' 
case, these modules  $M$  have the property that one of 
these two classes of modules (namely an orthogonal class) 
is closed under direct summands, while the other does not 
necessarily have this closure property.  More precisely, 
according to  Bazzoni's paper  \cite{Ba}, quoted at the 
beginning of this section,  n - tilting modules (or just tilting
modules, for short), admit the following definition. 

\subsection{A short definition (for the general case).}  
Given an $R$-module  $M$  and a natural number   $n > 0$, 
we denote by  $\Gen_n(M)$   the class of all modules   $X$                                
such that there is an exact sequence of the form
$$
	M(1) \To \dots \To M(n) \To X \To 0,
$$
where the  $M(i)$'s   are direct summands of direct sums of 
copies of   $M$.  Following  [Ba], we say that  $M$ is a  tilting  
(resp.~partial tilting) module of projective dimension at most  $n$     
if   $\Gen_n(M)$   is equal to (resp.~ is contained in) the orthogonal    
class    $M^\perp     =  \bigcap_{i>0} \ker  \Ext^i_R(M,-).$

\subsection{Example: A class of modules not always closed 
under direct summands.}  
Let  $R$  be the algebra given by 
the (fully commutative) quiver
$$
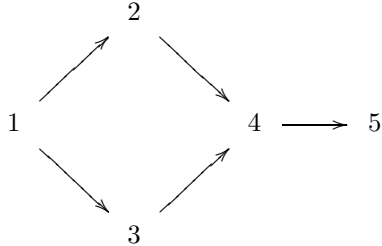
\xymatrix{
			&	\coluno{2} \ar[dr]	& & \\
\coluno{1} \ar[ur] \ar[dr] & 		& \coluno{4} \ar[r] & \coluno{5} \\
			&	\coluno{3} \ar[ur] &&  
}$$
such that the composition of any two arrows is zero.  Then we deduce from \cite[Example B]{D4}  that the injective module $ T  = \coldue{4}{5} \oplus \coldue{2\bx 3}{\bx4\bx} \oplus \coldue{1}{3} \oplus \coldue{1}{2}$ is a partial tilting module (of projective dimension  3 )  such that   $\Gen_3  (T)$
contains the module $1 \oplus 1$,  but not its summand  1.  
Moreover both  $\Gen_3(T)$   and  the class   $\Add (T)$, formed                                                         
by all injective  $R$--modules without simple summands, 
have the same indecomposable modules, namely the  4
indecomposable summands of  $T$.  This means that we 
cannot determine the class   $\Gen_3(T)$  by just looking at
the Auslander--Reiten quiver.  In other words, this means 
that the operation of making direct sums (that is the ``only 
really well - understood construction''  in the words of  \cite[page  476 ]{V})  is not enough to investigate an important class of modules generated by a maximal direct summand 
of a tilting module of projective dimension  3  .  

\subsection{Cancellation of an injective non projective 
summand (to obtain a ``large'' partial tilting module).}
We may roughly speaking say that we obtain the faithful 
module  $T$  constructed in  3.3  from the minimal injective  
cogenerator   $D(R_R)  = \Hom_K  (R_R  , K)$   by means of 
cancellation of its injective summand  1  (of projective                                                   
dimension  3).  On the other hand,  $\Ext^2_R(T , 5)\not = 0$ 
and   $5$  is the unique indecomposable module   $X$   such   
that   $\Hom_R  (T , X) = 0$.  Consequently, we have
\begin{claim}{$\star$}
 $\ker \Hom_R  ( T , - ) \  \bigcap  \    T ^\perp     =  0.$
\end{claim}
In other words,  $T$  is a large partial tilting module  \cite{D3}. 
We recall some facts concerning these modules.  First of 
all, any tilting module is a large partial tilting module \cite[page  371]{Ba}.  
Secondly, any finitely generated large partial 
tilting module of projective dimension at most  1  is a 
tilting module  \cite[Theorem 1]{C1} . Finally, property  ($\star$)  
implies that any large partial tilting module, say again  $T$, 
is sincere  \cite{KT} , that is with the property that  $\Hom (P,T) \not = 0$  
for every projective module  $P \not = 0$.  Consequently, a 
module  $T$  of finite length is sincere  (\cite{AuReS}  and  \cite{R4}) 
if every simple module is a composition factor of  $T$. 
After we point out the theoretical and practical importance 
of cancellations in Tilting Theory, we show that sometimes 
cancellation of a projective - injective summand (that is, of 
an obvious summand) of a tilting module gives rise to a 
large partial tilting module.

\subsection{Cancellation of a projective - injective summand 
         of a tilting module with minimal orthogonal class
   (to obtain a ``large'' partial tilting module). }
Let  $R$  be a finite dimensional algebra or, more generally, 
a noetherian and semiperfect ring such that every 
indecomposable injective module has a simple socle.  Let  
$M$  be an injective tilting $R$--module (of projective
dimension  $> 1$)  such that the orthogonal class
$M^\perp = \bigcap_{n>0} \ker\Ext^n_R(M,-)$ is the class of all injective
modules.  Let  $T$  be a sincere summand of  $M$, obtained 
from  $M$  after cancellation of a projective summand  $P$. 
Then we deduce from \cite[Theorem 4]{D3} that  $T$  is a large 
partial tilting module.  

We will give in  3.7  a ``minimal'' example of the above
result, where  $T$  is
\begin{enumerate}[(1)]
 \item the unique indecomposable injective module which 
is not projective;
\item a uniserial module such that every simple module has 
multiplicity one as a composition factor of  $T$ , that is a
sincere module of minimal dimension.
\end{enumerate}
  
\subsection{Remark on complexes.}  It turns out that various 
types of cancellations seem to be useful also by dealing 
with more abstract partial tilting objects.  For instance, 
this often happens with partial tilting complexes, say  $T^\circ$, 
in the sense of  Rickard \cite{Rk}  with the following property:
\begin{claim}{$a$}
 $T^\circ$  is the projective resolution of a large partial 
       modules, say  $T$ , which is not a tilting module. 
\end{claim}
By  \cite{Rk}, this hypothesis guarantees the existence of 
a non-zero right bounded complex (of projective modules)  $X^\circ$  
with the following property:
\begin{claim}{$b$}
 $X^\circ$   is not the projective resolution of a module
     and every morphism from  $T^\circ$  to any shift of  $X^\circ$
     is homotopic to zero. 
\end{claim}
The examples constructed in \cite{D5}, \cite{D6}  and \cite{D7} 
suggest that there is no canonical way to obtain  $X^\circ$  
from  $T^\circ$.  Moreover, the same holds by confining 
ourselves to complexes  $T^\circ$  and  $X^\circ$  satisfying ($a$)  
and  ($b$)  respectively and with the following
additional ``very combinatorial'' property (in the words 
of  \cite{Sc-ZI}) :
\begin{claim}{$c$}
Any non--zero component of the indecomposable 
summands of  $T^\circ$  and  $X^\circ$  is an indecomposable module.  
\end{claim}

Indeed, an indirect proof of the intricacy of complexes 
with respect to modules is that -- up to shift -- the choices 
of the indecomposable complexes  $X^\circ$  satisfying both ($b$)
and ($c$)  may be quite different.  For instance, they may 
be either zero  \cite[Example C  (iii)  and  (iv)]{D5}, or one  
\cite[Remark after Example 1]{D7}, or infinitely but countably 
many  \cite[Remark after Example 3]{D7} , or uncountably 
many  \cite[Example 4]{D7}.  The following example shows that 
by deleting some components of an indecomposable complex  
$T^\circ$, with properties  ($a$) and ($c$), we may obtain all the 
shortest complexes  $X^\circ$ with properties ($b$) and ($c$), that 
is the ``elementary'' complexes in the sense of  \cite{Sc-ZI} 
of the form  
\begin{claim}{$d$}
 $0\To P \To Q \To 0$ with $P$  and  $Q$ indecomposable projective modules.
\end{claim}
\subsection{ Example of cancellations concerning 
partial tilting modules and complexes.}
For every even integer $m>2$, there is a uniserial 
non faithful injective module   $T$  (of projective 
dimension $m$) such that we obtain all the 
indecomposable complexes  $X^\circ$  with the above 
properties  ($b$) and ($d$)  by means of various types 
of cancellations of some components of  $T^\circ$, that is 
left cancellations, right cancellations and, sometimes, 
also central cancellations.  
Indeed, let  $A$  be the Nakayama algebra, considered 
in  \cite{Ma2}, given by the quiver
$$\xymatrix{
\coltre{}{\bullet}{1} \ar[r]^{a_1} &
\coltre{}{\bullet}{2} \ar@{.}[r] &
\coltre{}{\bullet}{n-1}\ar[r]^{a_{n-1}} &
\coltre{}{\bullet}{n} \ar@/^2pc/[lll]^{a_n} 
}
$$
with relation    $a_n\cdot \dots \cdot a_1 = 0$,
 where   $2n = m + 2$.  Next, let  $T$ 
denote the injective module of the form  

$$\scalebox{0.7}{$\begin{array}{c}
 2\\
 3      \\
\vdots \\	
 n\\
 1       
\end{array}$.}$$


Then we obtain  $T$  from the minimal injective 
cogenerator  $D(A_A) = \Hom_K(A_A,K)$   after cancellation 
of its $n-1$ indecomposable projective summands. Thus 
we deduce from  3.5  (or from  \cite[Example 6]{D3}) that  
$T$ is a large partial tilting module of projective 
dimension  $m$.  Moreover its projective resolution  $T^\circ$  
satisfies  ($a$) and ($c$), and the complexes  $X^\circ$  satisfying  
($b$) and ($d$)  are of the form $0\To I(i)\To I(j) \To 0$,  
where  $I(*)$  denotes the indecomposable injective module corresponding to the vertex $*$ and  $i > j > 1$ \cite[Proposition 1 ]{D6}.  Hence they are exactly the complexes 
with two non - isomorphic injective components different 
from zero, obtained from  $T^\circ$  after suitable cancellations.    

\section{Do finite dimensional bimodules have a 
     concealed topology?}
Thanks to the method of visualizing bimodules by means
of pictures, as in Section 1, I could first ``see'' and then 
prove that even rather small bimodules have bad 
behaviour with respect to quite natural possible
constructions, namely embeddings into bimodules with 
an underlying left (or right) injective module.  Indeed, 
by the results proved in \cite{Ma1}  on the socle of  $E(C)/C$  
\cite[Lemma 2.2]{Ma1}  and on the modules cogenerated 
by  $E(C)/C$  \cite[Propositions  1.7  and  2.1  and 
Theorem 1.17]{Ma1}, it is natural to measure the gap 
between a classical cotilting module  $C$  and its 
injective envelope  $E(C)$, at least in case of modules 
which are finite dimensional vector spaces.  However, 
in this special situation, where the discrete topology 
should be the canonical topology, two radically 
different situations show up.

\subsection{Bad case.} It is not always possible to embed 
a finite dimensional cotilting bimodule  $C$  in another 
bimodule  $D$  with the property that   $D$, as a left (resp.~ 
right) module, is the injective envelope of  $C$ \cite[Example B (c), (d)]{D2}.  
Moreover, no left--right 
symmetry exists, because only one of the constructions 
may be possible \cite[Example A (c), (d)]{D2} .

\subsection{Good case.} When such an embedding exists, 
the structure of  $D$, namely its structure as a right 
(resp.~left)  module seems to be the most obvious one. 
Indeed also multiplications on the opposite side, that is 
right (resp.~left) multiplications, are described by nice 
matrices with many entries equal to zero. In other 
words, they seem to be ``continuous'' extensions of their 
restriction to  $C$.  However, as the following toy example 
shows, the property of being an indecomposable 
bimodule is neither hereditary nor left--right symmetric.

\subsection{Toy example of a finite dimensional cotilting 
          bimodule \cite[Example D]{D2}} 

Let $S$ (resp. $R$) be the algebra given by the quiver $\coltre{}{\bullet}{4}\To\coltre{}{\bullet}{5}\To\coltre{}{\bullet}{6}$ (resp. $\coltre{}{\bullet}{2}\To\coltre{}{\bullet}{1}\longleftarrow\coltre{}{\bullet}{3}$), and let
$_SC_R$ be an indecomposable cotilting bimodule with                                      
$$
_SC = \coltre{4}{5}{6} \oplus \coldue{5}{6} \oplus \coluno{5} \qquad 
\tonde{\text{resp.~} 
C_R = \coldue{1}{2\bx3} \oplus \coldue{1}{2} \oplus \coluno{2}
  }
 $$
 of the form

$$\text{\textbf{Picture 2}}\qquad \scalebox{0.3}{$
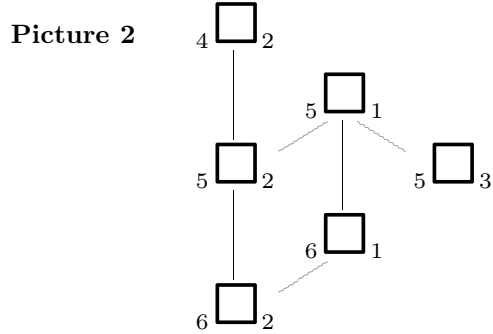
\xymatrix{
	 \quadr{4}{2} \ar@{-}[dd]  & &  \\
							    & \quadr{5}{1} \ar@{-}[dd]  \ar@{~}[dr] \ar@{~}[dl]   &  \\
	  \quadr{5}{2} \ar@{-}[dd] & & \quadr{5}{3}  \\
							    & \quadr{6}{1} \ar@{~}[dl]   &  \\
	 \quadr{6}{2}  &   
}$}
$$

Then  $C$  has codimension  2   in its left (resp.~ right) 
injective envelope 
$$E(C)  = \coltre{4}{5}{6} \oplus \coltre{4}{5}{6} \oplus \coldue{4}{5} 
\qquad 
\tonde{\text{resp.~} 
E(C) = \coldue{1}{2} \oplus \coldue{1}{2} \oplus\coldue{1}{2} \oplus\coldue{1}{3}}.
 $$ 
Moreover $E(C)$ is the support of an indecomposable (resp.~a decomposable) 
bimodule  $D$  , containing  $C$  as a bimodule, of the form

$$\text{\textbf{Picture 3}}\qquad \scalebox{0.3}{$
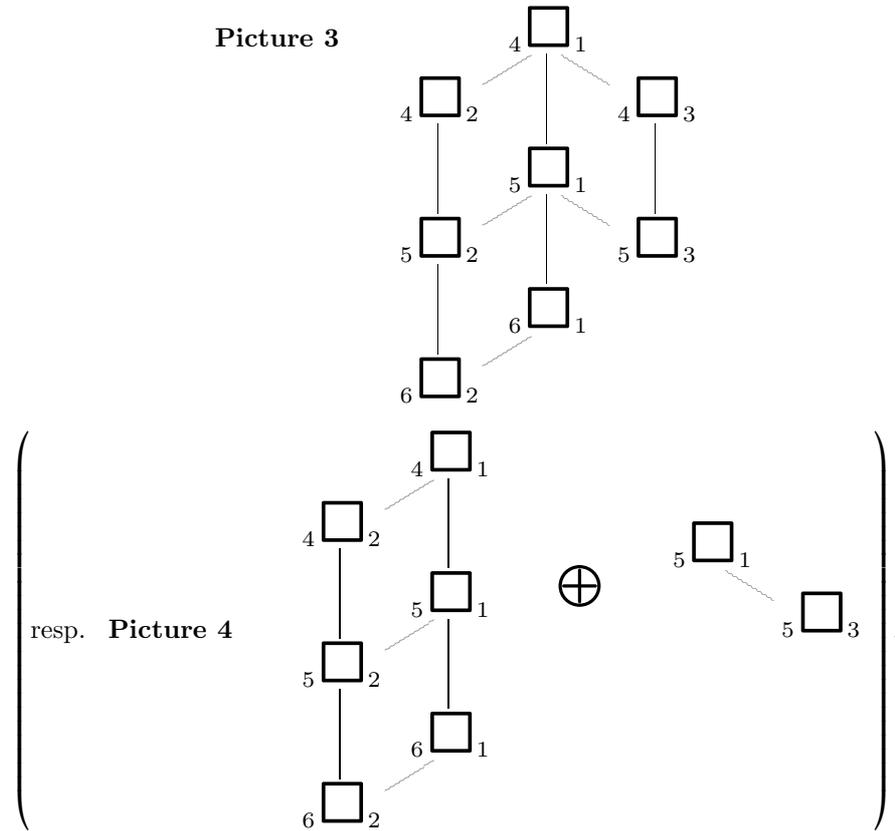
\xymatrix{
							    & \quadr{4}{1} \ar@{-}[dd]  \ar@{~}[dr] \ar@{~}[dl]   &  \\
	  \quadr{4}{2} \ar@{-}[dd] & & \quadr{4}{3} \ar@{-}[dd] \\
							    & \quadr{5}{1} \ar@{-}[dd]  \ar@{~}[dr] \ar@{~}[dl]   &  \\
	  \quadr{5}{2} \ar@{-}[dd] & & \quadr{5}{3}  \\
							    & \quadr{6}{1} \ar@{~}[dl]   &  \\
	 \quadr{6}{2}  &   
}$}
$$

$$
\left( \text{resp.~ \textbf{ Picture  4}}\qquad
\raisebox{22mm}{
\scalebox{0.3}{$
\xymatrix{
					& \quadr{4}{1} \ar@{-}[dd]  \ar@{~}[dl]     \\
	  \quadr{4}{2} \ar@{-}[dd] &  \\
							    & \quadr{5}{1} \ar@{-}[dd] \ar@{~}[dl]     \\
	  \quadr{5}{2} \ar@{-}[dd] & \\
							    & \quadr{6}{1} \ar@{~}[dl]     \\
	 \quadr{6}{2}  &   
}$}
}
\qquad 
\raisebox{5mm}{\scalebox{2.}{$\bigoplus$}}
\qquad
\raisebox{10mm}{
\scalebox{0.3}{$
\xymatrix{
	  \quadr{5}{1} \ar@{~}[dr] &  \\
							    & \quadr{5}{3}  }$} 
}\right).
$$

\subsection{Remarks  on the action of primitive idempotents 
         and nilpotent elements in the good case. } 
In all the examples constructed in  [D2] ,  where 
$_SC_R$ is a cotilting bimodule such that  $E(_SC)$  admits
a structure of  $S-R$  bimodule, containing the cotilting                                                                                     
bimodule     $_SC_R$    , the following facts hold:
\begin{enumerate}[(1)]
 \item The ring  $R$   is hereditary.
 \item If   $X$  is an indecomposable summand of  $_SC$    and  $e$                  
is a primitive idempotent of   $R$  such that   $xe = x$  for 
every element  $x\in X$, then we also have  $x' e = x'$     
for every element   $x' \in  E(_SX)$.
 \item The nilpotent elements of  $S$  and  $R$  corresponding 
to arrows  act on the elements of  $E(C)\setminus C$    in the 
easiest possible way, by a kind of shift.  Indeed, if   $s\in S$,  
$r\in R$  and  $u , v , w$   (resp.~  $u , w , w,  x$)  are linearly 
independent elements of  $C$  (resp.~ of  $E(C)$)  such that  
$x\in E(C)  \setminus  C$, $sx = v$, $sw = u$,  $vr = u$, then we 
have   $xr = w$, as illustrated by the following picture. 
\end{enumerate}

$$\text{\textbf{ Picture  5}}\qquad
\scalebox{0.3}{$
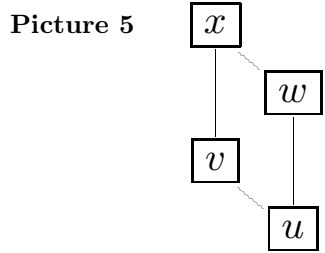
\xymatrix{
	  \scalebox{5.}{\fbox{$x$}} \ar@{-}[dd] \ar@{~}[dr] &  \\
							    & \scalebox{5.}{\fbox{$w$}}  \ar@{-}[dd]      \\
	  \scalebox{5.}{\fbox{$v$}}  \ar@{~}[dr]  & \\
							    & \scalebox{5.}{\fbox{$u$}}     }$
}
$$

\subsection{Remarks  on the action of new concealed rings 
in the bad case.}  In all the examples constructed in  \cite{D2},   
where $_SC_R$ is a cotilting bimodule such that  $E(_SC)$   does not admit 
a structure of  $S$--$R$  bimodule, containing the cotilting                                                               
bimodule   $_SC_R$, the following facts hold:
\begin{enumerate}[(1)]
 \item The ring  $R$   is not hereditary.
 \item There are a ring $R^\star$ , a ring epimorphism  $F : R^\star\To R$
and a bimodule   $_SU_{R^\star}$    , containing  $_SC_{R^\star}$, such that   $E(_SC) = _S U$.
 \end{enumerate}
We give two examples, where  $\ker F$  is a  $K$--vector space 
of dimension  1  or  2 .

\subsection{An example (of the bad case) with  $R^\star$ 
hereditary.}  As in \cite[Example A]{D2} , let  $R$  
(resp.~ $S$ )  be the algebra given by the quiver     
$$
\xymatrix{
	\coltre{}{\bullet}{1} \ar[r]^a & \coltre{}{\bullet}{2} \ar[r]^b & \coltre{}{\bullet}{3} 
} \text{with relation } ba=0 \qquad \tonde{\text{resp.~}
\xymatrix{
	\coltre{}{\bullet}{4} \ar[r]^a & \coltre{}{\bullet}{5} \ar[r]^b & \coltre{}{\bullet}{6}
}
}.
$$
                                                                                  
Let $_SC_R$ be the cotilting bimodule such that   $_SC =  \coluno{6} \oplus \coltre{4}{5}{6} \oplus \coluno{4}$ and $C_R    =   \coldue{3}{2}\oplus\coluno{2}\oplus\coldue{2}{1}$.  Next, let  $R^\star$   be the hereditary algebra given by the  Dynkin diagram 
$  \xymatrix{  \coltre{}{\bullet}{1} \ar[r]  & \coltre{}{\bullet}{2} \ar[r] & \coltre{}{\bullet}{3}}$.
Finally, let   $F : R^\star \To R$   be the obvious ring epimorphism.  
Then  $\dim \ker F = 1$   and  $C$  , regarded as a   $S$--$R^\star$  bimodule, 
is contained in the  $S$--$R^\star$  bimodule   $U$  ,  satisfying  (1)  and  (2),  described by the following picture  
                                                                         .                                                                        

$$\text{\textbf{Picture 6}}\qquad \scalebox{0.3}{$
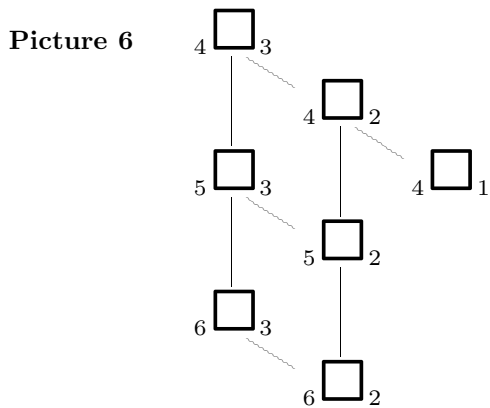
\xymatrix{
	  \quadr{4}{3} \ar@{-}[dd]  \ar@{~}[dr]   &  \\
							    & \quadr{4}{2} \ar@{-}[dd]  \ar@{~}[dr]   &  \\
	  \quadr{5}{3} \ar@{-}[dd] \ar@{~}[dr] & & \quadr{4}{1}  \\
							    & \quadr{5}{2} \ar@{-}[dd]    &  \\
	  \quadr{6}{3} \ar@{~}[dr] & &  \\
							    & \quadr{6}{2}    & 
}$}
$$

\bigskip

\subsection{An example (of the bad case)  with  R*   non
 hereditary.}  As in  \cite[Example C]{D2},  let  $R$  (resp.~ $S$ )  
be the  algebra given by the quiver 
$\xymatrix{
\coltre{}{\bullet}{1} \ar@/^/[r]^a &
\coltre{}{\bullet}{2} \ar@/^/[l]^b
} $ with relation $a b = 0$ (resp.~
$\xymatrix{
\coltre{}{\bullet}{3} \ar@/^/[r]^c &
\coltre{}{\bullet}{4} \ar@/^/[l]^d
} $ with relation $c d = 0$ ).  Next, let $_SC_R$ be the cotilting bimodule such that $_SC=\coltre{3}{4}{3}\oplus \coluno{3}$ and $C_R=\coltre{1}{2}{1} \oplus \coluno{1}$. Next, let $R^\star$ denote the algebra given by the quiver $\xymatrix{
\coltre{}{\bullet}{1} \ar@/^/[r]^a &
\coltre{}{\bullet}{2} \ar@/^/[l]^b
}$
with relation  $aba = 0$.  Finally, let $F : R^\star \To R$   be the 
obvious ring epimorphism.  Then  $\dim \ker F = 2$  and $C$, 
regarded as an $S$--$R^\star$ bimodule, is contained in the  
$S$--$R^\star$ bimodule $U$, satisfying  (1)  and  (2) , described 
by the following picture  

$$\text{\textbf{ Picture  7}}\qquad
\raisebox{22mm}{
\scalebox{0.3}{$
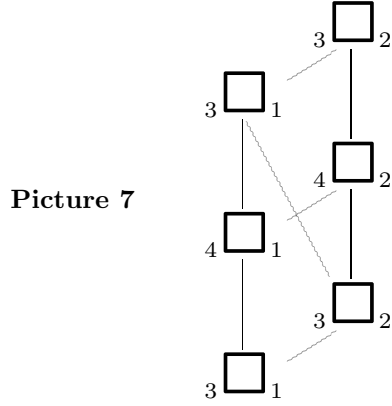
\xymatrix{
					& \quadr{3}{2} \ar@{-}[dd]  \ar@{~}[dl]     \\
	  \quadr{3}{1} \ar@{-}[dd] \ar@{~}[dddr] &  \\
							    & \quadr{4}{2} \ar@{-}[dd] \ar@{~}[dl]     \\
	  \quadr{4}{1} \ar@{-}[dd] & \\
							    & \quadr{3}{2} \ar@{~}[dl]     \\
	 \quadr{3}{1}  &   
}$}
}
$$

\bigskip

\subsection{ Two open problems on bimodules.}
  
\subsection*{Problem 1.}
Are conditions  (1) , (2)   and  (3)  of  Remark  
4.4	satisfied by any bimodule as in the good case described 
in  4.2?
\subsection*{Problem 2.}
 Are conditions  (1)  and (2)  of Remark 4.5
satisfied by any bimodule as in the bad case described in  4.1?

\section{Do infinite dimensional modules need no 
      topological tools?}
Another reason why the presence of a kind of topology,
pointed out in Section 4, is a strange fact is the absence
of topology in the proof of a result concerning dualities 
induced by cotilting bimodules of infinite dimension.
Before we discuss this, we recall that a left  $S$--module   
(resp.~a right  $R$--module)   $M$   is  reflexive with respect  
to  the bimodule    $_SU_R$    (or just  $U$--reflexive  or  reflexive, 
for short)   if    $M$  is canonically isomorphic to its
double dual with respect to  $U$ , that is to the group  
$ \Delta (\Delta(M))$,
 where  
$\Delta$ denotes both the contravariant functors $\Hom_?(-,_SU_R)$ for $?=R,S$ and the group $\Delta(X)$ 
%
is equipped with its bimodule structure 
(see \cite[Proposition 4.4]{AF}  or   
\cite[Propositions  3.4  and  3.5]{J}) for any left $S$-module and any right $R$-module $X$.

\subsection{Obvious and non obvious reflexive modules.}
Even in special cases, that is given a faithfully balanced 
bimodule $U$, there is a big gap between 
\begin{itemize}
\item[(*)] the well--known indecomposable reflexive modules, 
which are either projective or summands of  $U$  \cite[Propositions  20.13  and  20.14  and Corollary  20.16]{AF} ;
\item[(**)] the rest of the world, that is the non obvious 
indecomposable reflexive modules.
\end{itemize}   

Concerning  (*) , we first note that the cotilting bimodule 
described in  4.3  admits  4  indecomposable reflexive
left  (resp.~ right) modules.  Moreover, comparing 
Auslander--Reiten quivers (and looking at  Picture  2), 
we see that all of them are obvious reflexive modules and 
the duality $\Delta$  acts as follows:
\[\coltre{4}{5}{6}\mapsto \coluno{2}, \coldue{5}{6}\mapsto \coldue{1}{2\bx 3}, \coluno{5}\mapsto\coluno{3}, \coluno{6}\mapsto\coldue{1}{2}.\]
A similar situation holds for the cotilting module described in  
4.6 (resp.~ 4.7), where the cotilting duality $\Delta$ acts as follows:
\[\hspace*{-5mm}\coluno{6}\mapsto \coldue{3}{2}, \coltre{4}{5}{6}\mapsto \coldue{2}{1}, \coluno{4}\mapsto \coluno{1}, \coldue{5}{6}\mapsto \coluno{2}  
\tonde{\text{resp.~}
\coltre{3}{4}{3}\mapsto \coltre{1}{2}{1}, \coluno{3}\mapsto \coldue{2}{1}, \coldue{4}{3} \mapsto \coluno{1}\hspace*{-2mm}}.\]

However, to give an example of a cotilting bimodule 
admitting non obvious reflexive modules, it is enough
to take  Happel--Ringel's  cotilting (and tilting) bimodule     
$T$ described in Section 1.  Indeed, by comparing  
Auslander--Reiten quivers (and by looking at Picture 1), 
it is easy to see that $T$  admits   14  indecomposable 
reflexive left and right modules respectively.  On the
other hand,  $\coldue{5}{6}$  is the unique indecomposable projective
summand of the left  $A$--module  $T$.  Therefore $T$ 
admits   11   obvious (resp.~ 3  non obvious) 
indecomposable reflexive left and right modules 
respectively.  We already described (at the end of 
Section 1)  how the cotilting duality  $\Delta$  acts on the
indecomposable summands of $_AT$.  On the other hand,
$\Delta$  acts as follows on the remaining indecomposable 
obvious reflexive left  $A$--modules, that is on the five  
indecomposable projective modules which are not
summands  of $_AT$:    
$$\coltre{1}{2}{6}\mapsto \coldue{d}{e},\quad
\coldue{2}{6}\mapsto \coldue{b\bx{}d}{\bx{}e\bx},\quad
\coltre{3}{4}{6}\mapsto \coldue{d}{f}, \quad
\coldue{4}{6}\mapsto \coldue{d\bx{}c}{\bx{}f\bx},\quad
\coluno{6}\mapsto \coltre{\bx\bx\bx{}a\bx\bx\bx{}}{b\bx{}d\bx{}d\bx{}c}{\bx{}e\bx\bx\bx{}f\bx}.
$$
Finally, $\Delta$ acts as follows on the three indecomposable 
non obvious reflexive left  $A$--modules:
$$\coldue{2\bx5\bx4}{\bx6\bx6\bx} \mapsto \coldue{b\bx{}d\bx{}c}{\bx{}e\bx{}f\bx},\quad
\coltre{\bx\bx\bx\bx\bx3}{2\bx5\bx4}{\bx6\bx6\bx} \mapsto \coldue{b\bx{}d\bx{}}{\bx{}e\bx{}f}, \quad
\coltre{1\bx\bx\bx\bx\bx}{\bx2\bx5\bx4}{\bx\bx6\bx6\bx}\mapsto \coldue{\bx{}d\bx{}c}{e\bx{}f\bx}.
$$
Concerning  ($\star\star$)  at the beginning of this section, we refer to  
\cite[Sections 2 and 3 ]{C3} (or to \cite{C2,CbCF,CF} and to the 
other papers quoted in  \cite{C3})  for important results obtained 
by means of rather technical topological tools.  That's why it 
was a pleasant surprise to see that a discrete bimodule was 
enough to give the following answer to a question posed by 
Colpi with the hint:  ``You cannot use finite dimensional 
algebras and modules!''. (See  \cite[Theorem 1]{CbCF}  for 
the nice behaviour of submodules of reflexive modules over 
artin algebras.) 
                      
\subsection{Proposition \cite[Lemma 2.4  and  Theorem  2.5 (ii)]{D1}}       
Reflexive modules with respect to a cotilting bimodule are 
not necessarily closed under submodules.  Moreover, even 
the well--known reflexive modules with respect to a faithfully balanced bimodule $_SU_R$, that is the indecomposable summands of both the left (resp.~right)  regular module  $R$  
(resp.~ $S$)  and of the module   $U$   are not necessarily closed 
under submodules.

The next example - more precisely, the next picture - shows 
that the above result has a purely combinatorial motivation, 
coming from basic linear algebra.  Indeed, the result follows 
from the same reason why any infinite dimensional vector 
over  $K$  space cannot be isomorphic to its double dual, and 
so it cannot be reflexive with respect to the regular bimodule     
  $_KK_K$     (See  \cite[Proposition 1.8]{C3}  for a general result on                                        
direct sum of infinitely many non - zero reflexive modules 
with respect to a cotilting bimodule.)   
                                        
\subsection{Toy example of a finite dimensional cotilting 
bimodule  as in Proposition  5.2  \cite[Lemma 2.4  
and  Theorem  2.5  (ii)]{D1}.}  With terminology suggested by \cite{HU}, assume  $R = S$   is the 
``generalized'' Kronecker algebra given by the quiver
$$
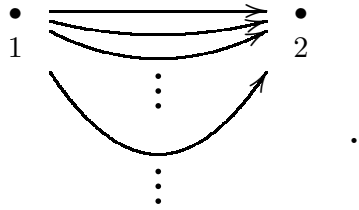
\scalebox{1.5}{
\xymatrix{
\coltre{}{\bullet}{1} \ar[rr]  \ar@/_/[rr] \ar@/_1pc/[rr]_{\medvdots}  
\ar@/_3pc/[rr]_\medvdots & &
\coltre{}{\bullet}{2}  
}\xymatrix{ & \\ \scalebox{1.2}{.}\\}}
$$
Then the indecomposable projective non simple left 
(resp.~ right) module is a reflexive module with respect 
to the cotilting bimodule  $_AA_A$.  However, as indicated                                                                            
in the following picture, its maximal submodule,  i.e.  the  
Jacobson radical of  $A$ , generated by the infinitely many 
arrows from  $\coluno{1}$  to   $\coluno{2}$,  is not reflexive with respect to $A$.

$$
\xymatrix{
	 \quadratino{1}{1} \ar[rd] \ar[rdd] \ar @<-1.5pc> @{} [rdd]^{\medvdots} & & 
	  \quadratino{2}{2} \ar[ld] \ar[ldd] \ar @<1.5pc> @{} [ldd]_{\medvdots} 
	  \\
	&   \quadratino{1}{2}  & \\
	&   \underset{\medvdots}{\quadratino{1}{2}}  
}
$$

\end{document}